\documentclass[runningheads]{llncs}
\usepackage{graphicx}
\usepackage{lmodern}
\usepackage[utf8]{inputenc}
\usepackage{fullpage}
\usepackage{amsmath}
\usepackage{hyperref}
\usepackage{csquotes}
\usepackage{microtype}
\usepackage{mathtools}
\usepackage{amssymb}
\usepackage{booktabs}
\usepackage{multirow}
\usepackage{wrapfig}
\usepackage[margin,noinline,draft]{fixme}
\usepackage{cryptocode}
\usepackage{tcolorbox}
\usepackage{tikz}
\usetikzlibrary{matrix}
\usetikzlibrary{graphs}
\usetikzlibrary{positioning, arrows.meta, calc}

\def\ve#1{\mathchoice{\mbox{\boldmath$\displaystyle\bf#1$}}
{\mbox{\boldmath$\textstyle\bf#1$}}
{\mbox{\boldmath$\scriptstyle\bf#1$}}
{\mbox{\boldmath$\scriptscriptstyle\bf#1$}}}

\newcommand\veb{{\ve b}}

\newcommand\vew{{\ve w}}
\newcommand\vex{{\ve x}}
\newcommand\vey{{\ve y}}

\newcommand\vezero{{\ve 0}}

\def\Z{\mathbb{Z}}
\def\N{\mathbb{N}}

\usepackage{stackengine}
\stackMath

\newcommand{\df}{:=}

\usetikzlibrary{fit,matrix,patterns,decorations.pathreplacing}

\newcommand{\Oh}{\mathcal{O}}
\newcommand{\OhOp}[1]{\Oh\mathopen{}\mathclose\bgroup\left( #1 \aftergroup\egroup\right)}

\usepackage{xspace}

\newcommand{\NPc}{\hbox{{\sf NP}-complete}\xspace}

\usepackage{tabularx}
\newcommand{\prob}[3]{
\begin{center}
\begin{tabularx}{\textwidth}{lX}
	\multicolumn{2}{l}{#1}\\
	{\bf Input:}&{#2}\\
	{\bf Find:}&{#3}
\end{tabularx}
\end{center}
}

\usepackage{xspace}

\newcommand{\II}{\mathcal{I}}

\newcommand{\KCCfree}{$(4K1, C_4, C_6)$-free\xspace}
\newcommand{\coloring}{\textsc{Graph Coloring}\xspace}
\newcommand{\nd}{\mathrm{nd}}

\AtEndEnvironment{proof}{\phantom{}\qed}

\begin{document}
\title{A Note on Coloring $(4K_1, C_4, C_6)$-free graphs with a $C_7$}

\authorrunning{Martin Koutecký}

\author{
Martin Kouteck{\'{y}}\inst{1}\thanks{Partially supported by Charles University project UNCE/SCI/004 and by the project 19-27871X of GA \v{C}R.}
}
\institute{Computer Science Institute, Charles University, Prague, Czech Republic\\
\email{koutecky@iuuk.mff.cuni.cz}
}

\maketitle

\begin{abstract}
Even-hole-free graphs are a graph class of much interest.
Foley et al.~[Graphs Comb. 36(1): 125-138 (2020)] have recently studied $(4K_1, C_4, C_6)$-free graphs, which form a subclass of even-hole-free graphs.
Specifically, Foley et al. have shown an algorithm for coloring these graphs via bounded clique-width if they contain a $C_7$.
In this note, we give a simpler and much faster algorithm via a more restrictive graph parameter, neighborhood diversity.
\end{abstract}

\section{Introduction}
A \textit{coloring} $c$ of a graph $G$ with $k$ colors is a mapping $c: V(G) \to \{1, 2, \dots, k\}$ such that no two neighboring vertices have the same color, and $\chi(G)$ is the smallest number such that $G$ admits a coloring with $\chi(G)$ colors.
The \coloring problem is as follows:

\prob{\coloring}{A graph $G$.}{A coloring of $G$ with $\chi(G)$ colors.}


This problem is \NPc in general, but its complexity is open for the well-studied~\cite{DBLP:journals/jct/ChangL15,DBLP:journals/jgt/ConfortiCKV02a,DBLP:journals/jgt/ConfortiCKV02,DBLP:journals/jct/SilvaV13} class of even-hole-free graphs, that is, graphs excluding all even cycles as induced subgraphs. 
Foley et al.~\cite{FoleyFHHL20} study \coloring on a subclass of even-hole-free graphs, namely graphs which do not contain a $4K_1$, $C_4$, nor $C_6$ as an induced subgraph, and they present some partial results.
One of them is showing that if such a graph does contain a $C_7$, then \coloring can be solved in polynomial time.
To show this, they argue that such graphs have cliquewidth at most $26$, and since \coloring can be solved in time $n^{f(k)}$ on graphs of cliquewidth $k$, this gives a polynomial algorithm for coloring these graphs.
It should be noted that the $f(k)$ dependence on $k$ is exponential and this cannot be improved under standard complexity assumptions~\cite{GolovachL0Z18}, so the complexity obtained by Foley et al.~\cite{FoleyFHHL20} could be as bad as $n^{2^{26}}$.

The purpose of this note is to point out that there is a much simpler and more efficient approach, which yields an algorithm of running time $\Oh(4^{26} \log n + n + m) = \Oh(n + m)$.
We build on the fact that Foley et al.~\cite{FoleyFHHL20} actually show much more than a cliquewidth bound: \cite[Theorem 2.6]{FoleyFHHL20} can be interpreted as saying that the neighborhood diversity of $G$ is at most $13$.
Neighborhood diversity is a graph parameter introduced by Lampis~\cite{Lampis}, who has also shown that \coloring can be solved in time $\Oh(f(\nd(G)) \log n + n + m)$ via integer programming.

In this note, we connect the notion of a $k$-uniform partition from~\cite{FoleyFHHL20} to neighborhood diversity, we describe the approach to solve \coloring, and we reduce the dependency on $\nd(G)$ from~\cite{Lampis} by using a more efficient algorithm for integer programming.

\section{Preliminaries}
For $n \in \N$, let $[n] = \{1, 2, \dots, n\}$.

\subsection{Neighborhood Diversity}
For a graph~$G$ we denote by $V(G)$ its set of vertices and $E(G)$ its set of edges.

Let $P_1, \dots, P_k$ be a partition of the vertex set $V(G)$ of a graph $G$.
We say that $(P_i)_{i \in [k]}$ is a \emph{$k$-uniform partition of $G$} if the subgraph of $G$ induced by $P_i$ is a clique for all $i \in [k]$, and for each $i,j \in [k], i \neq j$, between $P_i$ and $P_j$, there are either all possible edges or no edges in $G$.
If we change this definition slightly to allow $P_i$ to also be an independent set (not only a clique), we would get that $(P_i)_{i \in [k]}$ is a witness to the fact that the \emph{neighborhood diversity of $G$} is at most $k$, and we would call $(P_i)_{i \in [k]}$ an \emph{$\nd$-decomposition}.
To be precise, the neighborhood diversity of $G$, $\nd(G)$, is the smallest number $k$ such that an $\nd$-decomposition of $G$ of size $k$ exists.
It is also the number of equivalence classes of the relation $\sim \subseteq V(G)^2$ where $u \sim v$ iff $N(u) \setminus \{v\} = N(v) \setminus \{u\}$.
The neighborhood diversity of $G$ and a corresponding $\nd$-decomposition can be determined in linear time $\Oh(n+m)$~\cite{Lampis}.
Let $T(G)$ be a weighted multigraph which is defined according to the minimal $\nd$-decomposition of $G$ as follows.
The vertices of $T(G)$ are $v_1, \dots, v_k$, and the weight $n_i$ of $v_i$ is $|P_i|$. There is an edge $v_i v_j$, $i \neq j$, if the edges between $P_i$ and $P_j$ induce a complete bipartite graph in $G$.
There is a loop at $v_i$ if $P_i$ induces a clique in $G$.
We call $T(G)$ the \emph{type graph of $G$}.
Notice that the encoding length of $T(G)$ is $\Oh(k^2 \log n)$, and that it can be constructed in time $\Oh(n+m)$ from $G$.

\subsection{Integer Programming}
We write vectors in boldface (e.g., $\vex, \vey$) and their entries in normal font (e.g., the $i$-th entry of~$\vex$ is~$x_i$).
We use $\log \df \log_2$.
The \textsc{Integer Linear Programming} problem is to solve 
\begin{equation}
	\min \vew \vex:\, A\vex = \veb, \, \vex \geq \vezero, \, \vex \in \Z^{d},\label{IP} \tag{IP}
\end{equation}
where $\vew \in \Z^d$, $A \in \Z^{r \times d}$, and $\veb \in \Z^r$.
Let $\Delta = \|A\|_\infty = \max_{i,j} |A_{ij}|$ be the largest coefficient in $A$ in absolute value.
We will use a recent algorithm Jansen and Rohwedder~\cite{JansenR19} (which builds on~\cite{EisenbrandW20}).
\begin{proposition}[\cite{JansenR19}] \label{prop:JR}
\eqref{IP} can be solved in time $\Oh\left(\sqrt{r} \Delta\right)^{2r} \log \|\veb\|_\infty + \Oh(rd)$.
\end{proposition}

\section{Coloring \KCCfree Graphs with a $C_7$}
Lampis~\cite{Lampis} showed an algorithm solving \coloring in time $2^{2^{\Oh(\nd(G))}} \log n$ if $T(G)$ is provided.
Below, we present an improved approach to this problem, and then apply it to \KCCfree graphs with a $C_7$.

\begin{theorem} \label{thm:ndcoloring}
\coloring can be solved in time
$\Oh((\sqrt{\nd(G)})^{2\nd(G)} \log n)$ if $T(G)$ is given on input.
\end{theorem}
Note that the theorem above not only determines the correct value of $\chi(G)$ but also produces a valid coloring witnessing this fact.
\begin{proof}
Let $G$ be a graph with $\nd(G) = k$, $T(G)$ be its type graph, and recall that $n_i = |P_i|$ for each $i \in [k]$.
We define $\II(G) \subseteq 2^{V(T(G))}$ to be the set of maximal independent sets of $T(G)$, and we disregard the loops of $T(G)$ in this definition.
If $I'$ is a (not necessarily maximal) independent set in $T(G)$, we say that a $J \subseteq V(G)$ which is independent in $G$ \emph{induces} $I'$ in $T(G)$ if
\[
\left\{i \in [k] \mid J \cap P_i \neq \emptyset\right\} = \left\{i \in [k] \mid v_i \in I'\right\},
\]
and we say that $I$ \emph{dominates} $I'$ if $I' \subseteq I$.

We use the following key observation.
There exists an optimal coloring of $G$ such that each color intersects each clique $P_i$ at most once, and such that each $P_i$ inducing an independent set in $G$ is colored with just one color.
The first part is obvious; to see the second, simply choose any color appearing in $P_i$ in an optimal coloring and use it to color all of $P_i$.
We call any such coloring of $G$ \emph{canonical}.
Notice that each color class in a canonical coloring of $G$ induces an independent set $I'$ of $T(G)$ which is dominated by some maximal independent set $I \in \II(G)$.

Let us construct an~\eqref{IP} instance whose optimal value will be exactly $\chi(G)$ and whose optimal solution will encode an optimal canonical coloring of $G$.
There is a variable $x_I$ for each $I \in \II(G)$ whose intended meaning is how many colors in the optimal coloring induce some $I' \subseteq I$ of $T(G)$.
The constraints and objective function are then as follows:
\begin{align}
	\min & \sum_{I \in \II(G)} x_I & \label{eq:obj}\\
	\sum_{\substack{I \in \II(G)\\ i \in I}} x_I & \geq |P_i| & i \in [k] \text{ s.t. } G[P_i] \text{ is a clique} \label{eq:cover_clique} \\
	\sum_{\substack{I \in \II(G)\\ i \in I}} x_I & \geq 1 & i \in [k] \text{ s.t. } G[P_i] \text{ is an independent set} \label{eq:cover_indep} \\
	x_I & \geq 0 & \forall I \in \II(G) \label{eq:nonneg}
\end{align}
Consider an optimal solution of \eqref{eq:obj}--\eqref{eq:nonneg}.
We will now construct a canonical coloring of $G$ of value $OPT$, which is the optimal value of the objective function~\eqref{eq:obj}.
Initialize a counter $\chi := 1$.
Sequentially go over $I \in \II(G)$ and for each repeat $x_I$ many times the following:
\begin{enumerate}
	\item let $I' = \{i \in [k] \mid i \in I \wedge n_i \geq 1\}$, \label{it:1}
	\item let $J$ be an independent set of $G$ obtained by taking, for each $i \in I'$, one vertex of $P_i$ if $P_i$ is a clique, and all vertices of $P_i$ if $P_i$ is an independent set,
	\item color $J$ with color $\chi$ and remove $J$ from $G$,
	\item for each $i \in I'$, set $n_i := n_i - 1$ if $P_i$ is a clique and $n_i := 0$ if $P_i$ is an independent set,
	\item $\chi := \chi + 1$.
\end{enumerate}
Here, $I'$ is a maximal independent set of $T(G)$ dominated by $I$ such that there exists an independent set $J$ of $G$ which induces $I'$ in $T(G)$.
It is clear that our algorithm uses $OPT$ colors, and we have to argue that each vertex of $G$ is indeed colored.
To that end, view $\II(G)$ as a family of covering sets, $[k]$ as a universe of elements to be covered, and $n_i$ as a covering demand for element $i \in [k]$.
Interpret $x_I$ as the number of times we use $I \in \II(G)$ as a covering set.
The fact that constraints~\eqref{eq:cover_clique}--\eqref{eq:cover_indep} are satisfied means that each $i \in [k]$ is covered at least $n_i$ times.
This means that we could only ever attempt to color a vertex twice, but never that we would not color it at all.
We compensate for this potential ``overcovering'' by using an $I' \subseteq I$ instead of $I$ in step~\eqref{it:1} of the algorithm.

Conversely, we argue that each canonical coloring corresponds to some solution of~\eqref{eq:cover_clique}--\eqref{eq:nonneg}.
It is sufficient to follow the algorithm described above in reverse.
Initialize $\vex := \vezero \in \Z^{|\II(G)|}$.
Until $G$ is an empty graph, let $J$ be an independent set of $G$ corresponding to a color.
Let $I'$ be the independent set of $T(G)$ induced by $J$ and let $I \in \II(G)$ be arbitrary such that $I' \subseteq I$, increment $x_I := x_I + 1$, and remove $J$ from $G$.
Clearly,~\eqref{eq:obj} is the number of colors we have started with, and the constraints~\eqref{eq:cover_clique}--\eqref{eq:nonneg} are satisfied, so $\vex$ is a solution.

Finally, we solve the ILP~\eqref{eq:obj}--\eqref{eq:nonneg} using Proposition~\ref{prop:JR} which runs in time $(\sqrt{r} \Delta)^{2r} \log \|\veb\|_\infty + \Oh(rd)$ for an IP with $r$ rows, $d$ variables, largest coefficient $\Delta$, and with a right hand side $\veb$.
Here, $r = k$, $d = |\II(G)| \leq 2^k$, $\Delta = 1$, and $\log \|\veb\|_\infty \leq \log n$.
The resulting complexity is thus $(\sqrt{k})^{2k} \log n + \Oh(k 2^k) \leq \Oh((\sqrt{k})^{2k} \log n)$.

Note that the fact that constraints~\eqref{eq:cover_clique}--\eqref{eq:cover_indep} are inequalities while~\eqref{IP} requires equalities is no problem: we simply introduce $k$ so-called slack variables $s_1, \dots, s_k$, one per constraint.
For example, \eqref{eq:cover_clique} becomes $\sum_{\substack{I \in \II(G)\\ i \in I}} x_I + s_i = |P_i|$ and we have equality form, and this is at no asymptotical cost to our claims.
\end{proof}
Our approach improves on~\cite{Lampis} in two ways.
First, it replaces Lenstra's algorithm~\cite{Lenstra:1983} which depends exponentially on the number of variables $|\II(G)| \leq 2^{\nd(G)}$ by Proposition~\ref{prop:JR} which only depends exponentially on the number of rows $\nd(G)$.
Second, it only has a variable for each \emph{maximal} independent set, rather than for each independent set.
This can be helpful in situations where one has a lot of structural insight into $G$ and can explicitly bound the number of maximal independent sets $|\II(G)|$.

To apply Theorem~\ref{thm:ndcoloring} to \KCCfree graphs with a $C_7$, we need the fact that they have small $\nd(G)$. Foley et al.~\cite{FoleyFHHL20} show an even stronger result:
\begin{proposition}[{\cite[Theorem 2.6]{FoleyFHHL20}}]
	If $G$ is \KCCfree and contains a $C_7$, then it has a $k$-uniform partition with $7 \leq k \leq 13$.
\end{proposition}
Thus, we have these two straightforward corollaries:
\begin{corollary}
If $G$ is \KCCfree and contains a $C_7$, then $\nd(G) \leq 13$.
\end{corollary}
\begin{corollary} \label{cor:ndcoloring}
If $G$ is \KCCfree and contains a $C_7$, then \coloring can be solved on $G$ in time $\Oh(n + m)$, where the $\Oh$-notation hides a constant of order roughly $\sqrt{13}^{2 \cdot 13} \approx 2.9 \cdot 10^{14}$.
\end{corollary}
For contrast, the algorithm of Foley et al.~\cite{FoleyFHHL20} may have complexity as bad as $n^{2^{26}} \approx n^{6.7 \cdot 10^7}$.

\section{Discussion}
One may wonder whether the multiplicative constant in Corollary~\ref{cor:ndcoloring} could be improved.
The complexity of Proposition~\ref{prop:JR} is more precisely defined in terms of the so-called hereditary discrepancy of the constraint matrix. Because only relatively few type graphs $T(G)$ are possible for any \KCCfree graph $G$ containing a $C_7$, one could enumerate $\II(G)$ (which may be substantially smaller than the naïve bound $2^{\nd(G)}$), construct $A$, compute its hereditary discrepancy, and substitute this into the complexity statement.
Another approach to integer programming~\cite{EisenbrandHunkenschroederKleinKouteckyLevinOnn19} depends on the norm of so-called Graver basis elements; again, for a fixed matrix $A$, one can simply enumerate the Graver basis and check the norm.
Some experimental results~\cite{AltmanovaKK19} suggest that the worst-case upper bounds are often not attained, and much more practical constants are possible.

It is also worth noting that the approach described in Corollary~\ref{cor:ndcoloring} can be generalized to a wider graph class, namely graphs of bounded modular-width.
The definitions and full proof were given by Gajarský et al.~\cite{GajarskyLO13}, and we only give a brief outline here.
A graph has modular-width at most $k$ if it has a module decomposition with branching degree at most $k$; an inductive definition is that $G$ is either a clique or independent set, or it is possible to partition $G$ into sets $P_1, \dots, P_k$ such that between any $P_i, P_j$, $i \neq j$, there are either all possible edges or no edges, and each $G[P_i]$ has modular-width at most $k$.
The \coloring algorithm for graphs of small modular-width works in a bottom-up fashion along the module decomposition and replaces each (sub)module $G[P_i]$ with a clique $K_{\chi(G[P_i])}$.
When proceeding like this, each subproblem actually corresponds to coloring a graph with small neighborhood diversity, which can be done efficiently by Theorem~\ref{thm:ndcoloring}.

\section*{Acknowledgements}
Many thanks to Irena Penev for helpful discussions and feedback on a draft of this note.

\bibliography{milp}
\end{document}